\documentclass[12pt]{article}
\usepackage{amsmath,amssymb,amsthm,amscd,latexsym,}
 
\usepackage{graphicx}

\makeatletter
\renewcommand{\mod}[1]{\allowbreak \if@display \mkern 8mu \else
\mkern 5mu\fi {\operator@font mod}\,\,#1}
\makeatother

\newcommand{\bc}{\mathbb C}

\newcommand{\br}{\mathbb R}
\newcommand{\bz}{\mathbb Z}

\newcommand{\bp}{\mathbb P}

\newcommand{\Oc}{\mathcal O}

\newcommand{\aaa}{\mathbb A}
\newcommand{\ddd}{\mathbb D}
\newcommand{\eee}{\mathbb E}

\DeclareMathOperator{\Aut}{{\rm Aut}\,}

\DeclareMathOperator{\rk}{rk}

\newtheorem{theorem}{Theorem}

\begin{document}

\title{Some examples of K3 surfaces with infinite automorphism group which preserves an elliptic pencil}

\author{Viacheslav V. Nikulin}

\date{ }

\maketitle

\begin{abstract}
We give more details to our examples in \cite{Nik2} of K3 surfaces over $\bc$ such that
they have infinite automorphism group but it preserves some elliptic pencil of the K3.
\end{abstract}

\centerline{Dedicated to the memory of \`Ernest Borisovich Vinberg}

\section{Introduction}
\label{sec:introduction}
In this paper, we want to give more details to cases of \cite{Nik2} when there were
considered K3 surfaces $X$ which have the automorphism group $\Aut X$ which is infinite,
but preserves some elliptic pencil of the
K3 surface $X$. In Section \ref{sec:Nik2} we remind some results of \cite{Nik2} related to these cases.
In remaining sections we give more details to these cases considered in \cite{Nik2}.

\section{Some results of paragraph 4 of \cite{Nik2} when the automophism group of a K3 surface is infinite but it
preserves some elliptic pencil}
\label{sec:Nik2}

We remind to a reader (see \cite{Sh}, \cite{PS}) that a K3 surface is a K\"ahlerian compact complex surface $X$
such that its canonical class $K_X$ is zero
(i.e. there exists a holomorphic 2-dimensional differential form $\omega_X$ of $X$ with zero divisor) and $X$ has no non-zero  1-dimensional holomorphic differential forms (equivalently, $H^1(X, \Oc_X)=0$). In this paper, we shall consider
algebraic K3 surfaces.

In paragraph 4 of \cite{Nik2}, in particular, we gave examples of algebraic K3 surfaces such that their
automorphism group preserves some elliptic pencil of the K3 surface.

K3 surfaces $X$ of paragraph 4 of \cite{Nik2} have 2-elementary Picard lattices
$$
S=S_X =\{x\in H^2(X,\bz)\ |\ x\cdot \omega_X=0\}.
$$
Here a lattice means that $S$ is a free $\bz$-module of a finite rank
with an integer symmetric bilinear form which is given
by the intersection pairing in our case. The lattice $S$ is 2-elementary if
for $S^\ast=Hom(S, \bz)$ we have that the discriminant group
$$
A_S=S^\ast/S
$$
of the lattice $S$ is
2-elementary, that is $A_S\cong (\bz/2\bz)^a$ where $a\ge 0$ is an integer.

By global Torelli Theorem for K3 surfaces \cite{PS},
such K3-surfaces $X$ have and unique involution $\theta$ which
acts trivially (identically) on the Picard lattice $S$ but $\theta$ acts as $-1$ on the transcendental
lattice $T=T_X=(S_X)^\perp_{H^2(X,\bz)}$ of the K3 surface $X$. In particular, $\theta (\omega_X)=-\omega_X$ for a  2-dimensional holomorphic differential
form $\omega_X\in T_X\otimes \bc$ of $X$. Thus, $\theta$ is a canonical non-symplectic involution of $X$.
In particular, all automorphisms of $X$ commute with $\theta$.

Geometry of such K3 surfaces $X$ was studied in paragraph 4 of \cite{Nik2}, and their classification
was obtained.

For an algebraic K3 surface $X$, the Picard lattice $S=S_X$ is hyperbolic and has the signature
$(t_{(+)}=1, t_{(-)})$ where $t_{(+)}$, $t_{(-)}$ are the numbers of positive and negative squares respectively of
the real symmetric bilinear form $S\otimes \br$. Then $\rk S=r=1+t_{(-)}$. Moreover, $S$ is even, that is
$x^2=x\cdot x\equiv 0\mod 2$ is even for any $x\in S$.
There is the third important invariant $\delta$ of a 2-elementary lattice $S=S_X$ which is equal to $0$ or $1$.
The invariant $\delta$ is $0$ if $(x^\ast)^2\in \bz$ for any $x\in S^\ast$, otherwise, $\delta=1$.

The invariants $(r, a, \delta)$ define the isomorphism class of a 2-elementary Picard lattice $S=S_X$ of a K3 surface $X$.
See \cite{Nik1}

All possible $(r,a,\delta)$ are described in the table of the subsection 5 of the paragraph 4 of \cite{Nik2}:
we have
\begin{equation}
r+a\equiv 0\mod 2,\ \ r\ge 1, \ \ 0\le a\le r\le 20,\ \  1\le r+a\le 22,
\label{cond1}
\end{equation}
\begin{equation}
a\equiv 0\mod 2\ \  if\ \delta=0,\ \ r\equiv 2\mod 4\ \ if\ \ \delta=0,
\label{cond2}
\end{equation}
\begin{equation}
r\equiv 2\mod 8\ \ if\ \  a=0,\ \ r\equiv 2\pm 1\mod 8\ \  if\ \ a=1,
\label{cond3}
\end{equation}
\begin{equation}
\delta=1\ \ if\ \ r=a=6,\ \  \delta=1\ \ if\ \ (r,a)=(14,8).
\label{cond4}
\end{equation}

The set $X^\theta$ of fixed points of $\theta$ on $X$ was described.
The set $X^\theta$ is a non-singular algebraic curve of $X$
which may have several irreducible components which don't intersect each other.  If
$(r,a,\delta)\not= (10,8,0), (10,10,0)$, then
\begin{equation}
X^\theta=C^{(g)}+\sum_{i=1}^{k} E_i,\ \ g=11-(r+a)/2,\ \ k=(r-a)/2
\label{Xtheta1}
\end{equation}
where $C^{(g)}$ is a non-singular irreducible curve of genus $g$ and $E_i$, $1\le i\le k,$ are non-singular irreducible rational
curves (or $(-2)$-curves).
We have
\begin{equation}
X^\theta=C_1^{(1)}+C_2^{(1)}\ \ if \ \ (r,a,\delta)=(10,8,0)
\label{Xtheta2}
\end{equation}
where $C_1^{(1)}$ and $C_2^{(1)}$ are different non-singular
irreducible curves of genus $1$. We have
\begin{equation}
X^\theta=\emptyset\ \ if\ \ (r,a,\delta)=(10,10,0).
\label{Xtheta3}
\end{equation}

We remind to a reader that for a K3 surface $X$ for an irreducible curve $D$
we have
\begin{equation}
p_a(D)=\frac{D^2}{2}+1,
\label{pa}
\end{equation}
and
$E^2=-2$ and the linear system $|E|$ is zero-dimensional for a non-singular irreducible rational curve $E$,
and $C^2=0$ and the linear system $|C|$ is one-dimensional for a
non-singular irreducible curve $C$ of the genus $1$ (i.e., $C$ is an elliptic curve). Thus, $|C|$
defines an elliptic pencil.

Since the involution $\theta$ is canonical for $X$, it follows from \eqref{Xtheta1} --- \eqref{Xtheta3}
that the automorphism group $\Aut X$ preserves the elliptic pencil $\pi: X\to \bp^1=|C^{(1)}|$
(and the elliptic curve $C^{(1)}$)
in case \eqref{Xtheta1} if $r+a=20$, and the elliptic pencil $\pi:X\to \bp^1=|C_1^{(1)}|=|C_2^{(1)}|$
in case \eqref{Xtheta2} when $(r,a,\delta)=(10,8,0)$. These interesting cases when
the automorphism group $\Aut X$ of a K3 surface $X$ preserves an elliptic pencil $\pi: X\to \bp^1$ we want
to consider in more details in this paper. We remark that if $\Aut X$ is infinite then such invariant elliptic
pencil is unique. More generally, if an infinite subgroup $H\subset \Aut X$ preserves two different linear systems of elliptic curves $|C_1|$ and $|C_2|$, then $H$ is finite. Really, then $H$ preserves the
divisor class $h=cl(C_1+C_2)$ which has positive square $h^2$. From
the geometry of K3 surfaces, it follows that $H$ is finite.

In paragraph 4 of \cite{Nik2}, it was shown by considering degenerate fibers and the Mordell--Weyl group of the elliptic
pencil $\pi:X\to \bp^1$ that this pencil has no reducible fibers in the case \eqref{Xtheta1} when $k=0$,
equivalently, when $(r,a,\delta)=(10,10,1)$. In remaining cases,
the canonical involution $\theta$ has $2$ fixed points on $\bp^1$. In the case $\eqref{Xtheta2}$ when
$(r,a,\delta)=(10,8,0)$, the curves $C_1^{(1)}$ and $C_2^{(1)}$ give fibers of $\pi$ over these fixed points, and
in the case \eqref{Xtheta1} when $r+a=20$ and $(r,a,\delta)\not=(10,10,1)$, the fiber of $\pi$ over one fixed point is the elliptic curve $C^{(1)}$, and the fiber of $\pi$ over another
fixed point contains curves $E_i$, $1\le i\le k$, and is a reducible elliptic curve $\mathcal{E}$ of Dynkin type
\begin{equation}
\mathcal{E}=
\left\{
\begin{array}{ll}
\widetilde{\eee}_6 & \rm{if}\ k=4,\ \delta=0\\
\widetilde{\aaa}_{2k-1} &  \rm{in\ other\ cases}
\end{array}
\right.\ \ .
\label{Xtheta4}
\end{equation}
We recall that $k=(r-a)/2$.

In the next sections we shall give more details in each of these cases $(r,a,\delta)$
related to the elliptic curves. In particular, we calculate exactly classes in Picard lattice of curves $C^{(1)}$,
$C^{(1)}_1$, $C^{(1)}_2$ and $E_1,\dots,  E_k$
in \eqref{Xtheta1}, \eqref{Xtheta2}. We will give more details to the proof of \eqref{Xtheta4} and
we calculate exactly classes in Picard lattice of components of $\mathcal{E}$ in \eqref{Xtheta4}.

\section{The case $(r,a,\delta)=(10,10,1)$}
\label{sec:10,10,1}
We shall use standard notations for lattices. By $\langle A \rangle$,
we denote a lattice with the integer symmetric matrix $A$. By $M(n)$
we denote the lattice which is obtained from a lattice
$M$ by multiplication of the symmetric bilinear form of $M$ by the interger $n$.
By $nM$ we denote the orthogonal sum of $n$ copies of a lattice $M$ where $n\ge 0$ is
an integer.
By $M_1\oplus M_2$ we denote the orthogonal sum of
lattices $M_1$ and $M_2$. By $A_n$, $n\ge 1$, $D_m$, $m\ge 4$, $E_k$, $k=6,7,8$,
we denote standard negative definite root lattices
correspoinding to Dynkin diagrams $\aaa_n$, $\ddd_m$, $\eee_k$ respectively with roots with square $(-2)$. See \cite{Bou}.

The even hyperbolic $2$-elementary lattice with the invariants $(r=10,a=10,\delta=1)$ is
isomorphic to the lattice
$$
S_{10,10,1}=
\langle\left(\begin{array}{rr}
 0 & 2 \\
 2 & -2
 \end{array}\right)\rangle\oplus E_8(2)
$$
with the bases $\{c,d,e_1,e_2,e_3,e_4,e_5,e_6,e_7,e_8\}$
where $\{c,d\}$ have the matrix
$\left(\begin{array}{rr}
 0 & 2 \\
 2 & -2
 \end{array}\right)$ and $e_1,...,e_8$ give the standard basis
 of $E_8$ with the diagram $\eee_8$ (see \cite{Bou}).

 We consider K3 surfaces $X$ with the Picard lattice $S_{10,10,1}$,
 elliptic pencil $|C|$ with the $cl(C)=c$ and non-singular irreducible rational curve $D$ with the class $d$.
 Such $X$ exist,  since $E_8(2)$ has no elements with square $-2$ and $c$, $d$ define
 the NEF cone of $X$.

 Since $(c)^\perp_S=\bz c\oplus E_8(2)$ has no elements with square $-2$, it
 follows that the Mordell---Weil group of this elliptic pencil $|c|$ (see \cite{Sh})
 contains the subgroup $\bz^8$, $8=\rk E_8(2)$, which gives a subgroup of finite index of
 the automorphism group of the elliptic pencil $|C|$. Since this group is infinite,
 $c=cl(C^{(1)})$ and $k=0$ in \eqref{Xtheta1}. The elliptic pencil $|c|$ is the unique
 elliptic pencil of $X$ with infinite automorphism group and $\Aut X$ is the automorphism
 group of this elliptic pencil. We remark that $C^{(1)}\cap D$ give two fixed points of
 the canonical involution $\theta$ on the non-singular rational curve $D$.

\section{The case $(r,a,\delta)=(10,8,0)$}
\label{sec:10,8,0}

An even hyperbolic 2-elementary lattice with the invariants $(r=10,a=8,\delta=0)$
is isomorphic to
$$
S_{10,8,0}=U\oplus E_8(2)\ where \
U=\langle\left(\begin{array}{rr}
 0 & 1 \\
 1 & -2
 \end{array}\right)\rangle\
$$
with the bases $\{c,d,e_1,e_2,e_3,e_4,e_5,e_6,e_7,e_8\}$
where $\{c,d\}$ have the matrix $U$
and $e_1,...,e_8$ correspond to the standard basis
 of $E_8$ with the diagram $\eee_8$ (see \cite{Bou}).

We consider K3 surfaces $X$ with the Picard lattice $S_{10,8,0}$,
elliptic fibration $|C|$ with the $cl(C)=c$ and non-singular irreducible rational curve $D$ with the class $d$.
We write fibration instead of pencil because for this case $D$ gives a section
of this fibration: we have $c\cdot d=1$ and the linear system $|c|$ can be identified with $D$ .

 Such $X$ exist,  since $E_8(2)$ has no elements with square $-2$ and $c$, $d$ define
 the NEF cone of $X$.

Since $(c)^\perp_S=\bz c\oplus E_8(2)$ has no elements with square $-2$, it
 follows that the Mordell---Weil group of this elliptic fibration $|c|$ (see \cite{Sh})
 contains the subgroup $\bz^8$, $8=\rk E_8(2)$, which gives a subgroup of finite index of
 the automorphism group of the elliptic fibration $|C|$. Since this group is infinite,
 $c=cl(C^{(1)})=cl(C^{(2)})$ in \eqref{Xtheta2}. The elliptic fibration $|c|$ is the unique
 elliptic pencil of $X$ with infinite automorphism group and $\Aut X$ is the automorphism
 group of this elliptic fibration. We remark that $C^{(1)}\cap D$ and $C^{(2)}\cap D$
 give two fixed points of the canonical involution $\theta$ on
 the non-singular rational curve $D$.

 This lattice $S_{10,8,0}=U\oplus E_8(2)$ was missed in Theorem 5.12 of the preprint \cite{Me}.

\section{The case $(r,a,\delta)=(11,9,1)$}
\label{sec:10,8,0}

An even hyperbolic 2-elementary lattice with the invariants $(r=11,a=9,\delta=1)$
is isomorphic to
$$
S_{11,9,1}=U\oplus E_8(2)\oplus A_1\ where \
U=\langle\left(\begin{array}{rr}
 0 & 1 \\
 1 & -2
 \end{array}\right)\rangle\
$$
with the bases $\{c,d,e_1,e_2,e_3,e_4,e_5,e_6,e_7,e_8,f_1\}$
where $\{c,d\}$ have the matrix $U$, elements
$e_1,...,e_8$ correspond to the standard basis
 of $E_8$ with the diagram $\eee_8$ (see \cite{Bou})
 and $f_1$ gives the basis of $A_1$, $f_1^2=-2$ .

We consider K3 surfaces $X$ with the Picard lattice $S_{11,9,1}$,
elliptic fibration $|C|$ with the $cl(C)=c$ and non-singular irreducible rational curve $D$ with the class $d$
which gives the section of this fibration, and non-singular irreducible rational curve $F_1$ with the class $f_1$.
Non-singular irreducible rational curves with classes $f_1$ and $f_2=c-f_1$ define a reducible fiber of
the type $\widetilde{\aaa}_1$ of this fibration with non-singular irreducible rational curves $F_1$ with $cl(F_1)=f_1$
and $F_2$ with $cl(F_2)=f_2=c-f_1$.

 Such $X$ exist,  since $E_8(2)$ has no elements with square $-2$ and $c$, $d$, $f_1$, $f_2=c-f_1$ define
 the NEF cone of $X$.

 Since $(c, f_1, f_2)^\perp_S=\bz c\oplus E_8(2)$ has no elements with square $-2$, it
 follows that the Mordell---Weil group of this elliptic fibration $|c|$ (see \cite{Sh})
 contains the subgroup $\bz^8$, $8=\rk E_8(2)$ which gives a subgroup of finite index of
 the automorphism group of the elliptic fibration $|C|$. Since this group is infinite,
 $c=cl(C^{(1)})$ in \eqref{Xtheta1}. The elliptic fibration $|c|$ is the unique
 elliptic pencil of $X$ with infinite automorphism group and $\Aut X$ is the automorphism
 group of this elliptic fibration. We remark that $C^{(1)}\cap D$ and $F_2\cap D$
 give two fixed points of the canonical involution $\theta$ on
 the non-singular rational curve $D$. In \eqref{Xtheta1} we have $k=1$ and
 $E_1$ has the class $f_2=c-f_1$ and $E_1=F_2$ . The curve $\mathcal{E}$ in \eqref{Xtheta4}
 has the type $\widetilde{\aaa}_1$ with components
 with classes $f_1$ and $f_2=c-f_1$.


\section{Cases $(r,a=20-r,\delta=1)$,  $r=12,...,17$,\\ and $(r=18,a=2,\delta=0)$}
\label{sec:generalcase}

These are cases when
$$
S=U\oplus D_{16-2t}\oplus tA_1,\ \ \ t=0,1,2,3,4,5,6.
\label{latS:caser.a.1}
$$
Then $r=18-t$, $a=2+t$ and $\delta=1$ if $t>0$, and $\delta=0$ if $t=0$.

The lattice $S$ has the basis $\{e,d,f_1,f_2,f_3,\dots,f_{16-2t},g_1,\dots,g_t\}$
where $\{e,d\}$ give the standard basis of $U$,  $f_1,f_2,f_3,\dots,f_{16-2t}$ is the
standard basis for the root lattice $D_{16-2t}$
and $\{g_1,\dots,g_t\}$ give the standard basis for the root lattice $tA_1$. See \cite{Bou}.

We consider K3 surfaces $X$ with the Picard lattice $S$, the elliptic fibration $|e|$ with the section
$D$ which is a non-singular irreducible rational curve with the class $d$, the degenerate
fiber of the type $\widetilde{\ddd}_{16-2t}$ with the divisor $F_0+F_1+2F_2+2F_3+\cdots +2F_{14-2t}+F_{15-2t}+F_{16-2t}$ where non-singular irreducible rational curves $F_i$, $i=1,...,{16-2t}$, have classes $f_i$ respectively and $F_0$ has the class $f_0=e-f_1-2f_2-2f_3-2f_4-\cdots - 2f_{14-2t}-f_{15-2t}-f_{16-2t}$, and $t$
degenerate fibers of the type $\widetilde{\aaa}_1$ with the divisors $G_j+G_j^\prime$, $j=1,\dots,t$,
 where $G_j$ has the class $g_j$, and $G_j^\prime$ has the class $e-g_j$. Curves $G_j$ and $G_j^\prime$ are also
 non-singular irreducible rational curves.

The canonical involution $\theta$ preserves all non-singular irreducible rational curves (with square $-2$) of $X$
 because they are defined by their classes in Picard lattice. We say that such a curve has the type $+$ if $\theta$ is
 identical on this curve (or it is one of the curves $E_1,\dots ,E_k$ in \eqref{Xtheta1}), and it has the type $-$ if not. Then $\theta$ has two fixed points on this curve. If a non-singular irreducible rational curve of the type $+$ intersects  transversally another non-singular irreducible rational curve, then the second curve has the type $-$ becase of the action
 of $\theta$ by $(-1)$ on the canonical form $\omega_X$. It follows that the curves $D$ with the class $d$ and
 $F_2, F_4, F_6,\dots,F_{14-2t}$ with the classes
 $f_2, f_4, f_6,\dots, f_{14-2t}$ respectively have the type $+$. Thus they give the curves $E_1,E_2,\dots, E_k$ where $k=(r-a)/2=8-t$ in \eqref{Xtheta1}. Thus, they have classes
 \begin{equation}
 d, f_2, f_4, f_6,\dots, f_{14-2t}.
 \label{+:general}
 \end{equation}
 All other non-singular irreducible rational curves on $X$ have the type $-$. Let $C=C^{(1)}$ in
 \eqref{Xtheta1}. Two fixed points of $\theta$ on a non-singular irreducible rational curve $M$ on $X$
of the type $-$
 come from intersections of $M$ with $D$, $F_2, F_4, \dots , F_{14-2t}$ and $C$.
 Using this information, we obtain that the class $c$ of $C$ is
 equal to
$$
c=3e^\ast+f_1^\ast+f_{15-2t}^\ast+f_{16-2t}^\ast+2g_1^\ast+2g_2^\ast+ \cdots +2g_t^\ast=
$$
\begin{equation}
6e+3d-2f_1-3f_2-4f_3-5f_4-6f_5-\cdots -(15-2t)f_{14-2t}-(8-t)f_{15-2t}-(8-t)f_{16-2t}-g_1-\cdots -g_t
\label{c:general}
\end{equation}
where we denote by
$\{e^\ast,d^\ast,f_1^\ast,\dots ,  f_{16-2t}^\ast, g_1^\ast,\dots, g_t^\ast\}$
the corresponding dual basis from $S^\ast$ for our  basis in the lattice $S$. One can check that $c^2=0$.

Our considerations show that the curve $D$ is the unique section of the elliptic fibration $|e|$ because we
showed that this section is one of $k$ curves $E_1,\dots , E_k$ from \eqref{Xtheta1}, but other of these $k$ curves
are classes $f_2, f_4,\dots, f_{14-2t}$ of the degenerate fiber of the type $\widetilde{\ddd}_{16-2t}$ of $|e|$ which are distinguished in the Dynkin diagram $\widetilde{\ddd}_{16-2t}$.

The elliptic pencil $|c|$ has the section $F_1$ because $f_1\cdot c=1$. Since the involution $\theta$
has exactly two fixed points in $F_1$, one of this points belongs to the fiber $C$ and another fixed point belongs to
the fiber of $|c|$ with
curves
$D, F_0, F_2, F_3,\dots, F_{14-2t}$ with classes $d, f_0, f_2, f_3,\dots , f_{14-2t}$ respectively.
The class
$$
\alpha=c-d-f_0-f_2-\dots - f_{14-t}
$$
has $\alpha^2=-2$ and it has non-negative intersection with all
classes of non-singular irreducible rational curves $M$ of $X$ which have $e\cdot M=0$ or $e\cdot M=1$.
We have $e\cdot \alpha=2$. It follows that $\alpha$ corresponds to an irreducible non-singular rational curve on $X$
by Vinberg's criterion \cite{Vin}. It follows that
\begin{equation}
d, f_0, f_2, f_3,\dots , f_{14-2t}, \alpha
\label{calE:general}
\end{equation}
give all classes of irreducible components of the degenerate fiber of the type $\tilde{\aaa}_{15-2t}$ of the elliptic fibration $|c|$. The elliptic fibration $|c|$ has no other
reducible fibers because they contain non-singular rational curves and should be over other fixed points of $\theta$
on the section $F_1$ of $|c|$. But $\theta$ has only $2$ fixed points on $F_1$.

Thus, from the theory of elliptic pencils on K3 we obtain that the subgroup of finite index of the automorphism group of
the elliptic fibration $|c|$ is a subgroup of finite index of the group
$MW^\ast$ where
$$
MW=(c,f_1,f_3,f_4, \dots f_{14-2t},f_0,d,\alpha)^\perp_S.
$$
which has the rank $t+1=a-1$. The lattice $MW$ has no elements with square $(-2)$.

Thus, we obtain

\begin{theorem}
The automorphism group $\Aut X$ of a K3 surface over $\bc$ with the 2-elementary
Picard lattice $S$ with the invariants
$(r,a=20-r,\delta=1)$ where $r=12,...,17$, and $(r=18,a=2,\delta=0)$, equivalently,
$S=U\oplus D_{16-2t}\oplus tA_1$ where $t=a-2$,
is equal to
the automorphism group of the elliptic fibration $|c|$ of $X$ where $c$ is given in  \eqref{c:general}.
Classes of the curves $E_1,\dots E_k$ from \eqref{Xtheta1}
are given in \eqref{+:general} where $k=(r-a)/2=8-t$.
The elliptic fibration $|c|$ has
the degenerate fiber of the type $\widetilde{\aaa}_{15-2t}$ with classes \eqref{calE:general} of
components and no other reducible fibers. Then $\Aut X$
is isomorphic to $\bz^{a-1}=\bz^{t+1}$ up to finite index.

See details in the considerations above.
\label{theorem:general}
\end{theorem}

\section{The case $(r,a,\delta)=(18,2,1)$}
\label{sec:18,2,1}

In this case,
\begin{equation}
S=S_{18,2,1}=U\oplus E_8\oplus E_7\oplus A_1.
\label{latS:18.2.1}
\end{equation}

The lattice $S$ has the basis $\{e,d,f_1,f_2,f_3,f_4,f_5,f_6,f_7,f_8,g_1,g_2,g_3,g_4,g_5,g_6,g_7,h_1\}$
where $\{e,d\}$ is the standard basis of $U$,  $f_1,\dots,f_8$ is the
standard basis of $E_8$, $g_1,\dots,g_7$ is the standard basis of $E_7$ and $h_1$ is the standard basis of $A_1$.
See \cite{Bou}.

We consider K3 surfaces $X$ with the Picard lattice $S$, the elliptic fibration $|e|$ with the section
$D$ which is a non-singular irreducible rational curve with the class $d$, the degenerate
fiber of the type $\widetilde{\eee}_8$ with the divisor $2F_1+3F_2+4F_3+6F_4+5F_5+4F_6+3F_7+2F_8+F_0$ where non-singular irreducible rational curves $F_i$, $i=1,...,8$, have classes $f_i$ respectively and $F_0$ has the class $f_0=e-2f_1-3f_2-4f_3-6f_4-5f_5-4f_6-3f_7-2f_8$, degenerate fiber of the type $\widetilde{\eee}_7$ with the
divisor $G_0+2G_1+2G_2+3G_3+4G_4+3G_5+2G_6+G_7$ where $G_i$, $i=1,\dots 7$, give the standard basis of $E_7$
and $G_0$ has the class $g_0=e-2g_1-2g_2-3g_3-4g_4-3g_5-2g_6-g_7$, degenerate fiber of the type
$\widetilde{\aaa}_1$ with the divisors $H_1+H_1^\prime$,
 where $H_1$ has the class $h_1$, and $h_1^\prime$ has the class $e-h_1$. Curves $H_1$ and $H_1^\prime$ are also
 non-singular irreducible rational curves.

By considering the canonical involution $\theta$, like for the previous case, we obtain that
the curves $F_1, F_4, F_6,F_8, D, G_1, G_4, G_6$  with the classes
\begin{equation}
f_1, f_4, f_6,f_8, d, g_1, g_4, g_6
\label{+:18.2.1}
\end{equation}
respectively
have the type $+$. Thus they give the curves $E_1,\dots E_k$ where $k=(r-a)/2=8$ in \eqref{Xtheta1}.
 All other non-singular irreducible rational curves on $X$ have the type $-$. Let $C=C^{(1)}$ in
 \eqref{Xtheta1}. Like for the previous case,
 using this information, we obtain that the class $c$ of $C$ is
 equal to
$$
c=3e^\ast+f_2^\ast+g_2^\ast+g_7^\ast+2h_1^\ast=
$$
\begin{equation}
6e+3d-5f_1-8f_2-10f_3-15f_4-12f_5-9f_6-6f_7-3f_8-3g_1-5g_2-6g_3-9g_4-7g_5-5g_6-3g_7-h_1
\label{c:18.2.1}
\end{equation}
where we denote by
$\{e^\ast,d^\ast,f_1^\ast,\dots,f_8^\ast, g_1^\ast,\dots, g_7^\ast, h_1^\ast\}$
the corresponding dual basis from $S^\ast$ for our  basis in the lattice $S$. One can check that $c^2=0$.

Our considerations show that the curve $D$ is the unique section of the elliptic fibration $|e|$ because we
showed that this section is one of $k$ curves $E_1,\dots , E_k$, $k=8$, from \eqref{Xtheta1}, but other of these $8$ curves
have classes $f_1, f_4,f_6,f_8,g_1,g_4, g_6$ of the degenerate fibers of types $\widetilde{\eee}_8$ and
$\widetilde{\eee}_7$ of $|e|$ which are distinguished in their Dynkin diagrams.

The elliptic pencil $|c|$ has the section $G_7$ because $g_7\cdot c=1$. Since the involution $\theta$
has exactly two fixed points in $G_7$, one of this points belongs to the fiber $C$ and another fixed point belongs to
the fiber of $|c|$ with
curves $D, F_1,F_3,F_4,F_5,F_6,F_7,F_8,F_0,G_0,G_1,G_3,G_4,G_5,G_6$ with classes
$d, f_1, f_3,f_4,f_5,f_6,f_7,f_8,f_0,g_0,g_1,g_3,g_4,g_5,g_6$ respectively.
The class
$$
\alpha=c-d-f_1-f_3-f_4-f_5-f_6-f_7-f_8-f_0-g_0-g_1-g_3-g_4-g_5-g_6
$$
has $\alpha^2=-2$ and it has non-negative intersection with all
classes of non-singular irreducible rational curves $M$ of $X$ which have $e\cdot M=0$ or $e\cdot M=1$.
We have $e\cdot \alpha=2$. It follows that $\alpha$ corresponds to an irreducible non-singular rational curve on $X$
by Vinberg's criterion \cite{Vin}. It follows that
\begin{equation}
d, f_1, f_3,f_4,f_5,f_6,f_7,f_8,f_0,g_0,g_1,g_3,g_4,g_5,g_6,\alpha
\label{Ecal:18.2.1}
\end{equation}
give classes of all irreducible components of the degenerate fiber of the type $\widetilde{\aaa}_{15}$ of the elliptic fibration $|c|$. The elliptic fibration $|c|$ has no other
reducible fibers because they contain non-singular rational curves and should be over other fixed points of $\theta$
on the section $G_7$ of $|c|$. But $\theta$ has only $2$ fixed points on $G_7$.

Thus, from the theory of elliptic pencils on K3 we obtain that the subgroup of finite index of the automorphism group of
the elliptic fibration $|c|$ is a subgroup of finite index of the group
$MW^\ast$ where
$$
MW=(c,g_7,d, f_1, f_3,f_4,f_5,f_6,f_7,f_8,f_0,g_0,g_1,g_3,g_4,g_5,\alpha)^\perp_S.
$$
which has the rank $1$. The lattice $MW$ has no elements with square $(-2)$.

Thus, we obtain

\begin{theorem}
The automorphism group $\Aut X$ of a K3 surface $X$ over $\bc$ with the 2-elementary
Picard lattice $S$ with the invariants
$(r=18,a=2,\delta=1)$, equivalently,
$S=U\oplus E_8\oplus E_7\oplus A_1$,
is equal to the automorphism group of the elliptic fibration $|c|$ of $X$
where $c$ is given in  \eqref{c:18.2.1}.
Classes of the curves $E_1,\dots E_k$ from \eqref{Xtheta1}
are given in \eqref{+:18.2.1} where $k=8$.
The elliptic fibration $|c|$ has
the degenerate fiber of the type $\widetilde{\aaa}_{15}$
with components with classes \eqref{Ecal:18.2.1}
and no other reducible fibers. Then $\Aut X$
is isomorphic to $\bz$ up to finite index.

See details in the considerations above.
\label{theorem18.2.1}
\end{theorem}


\section{The case $(r,a,\delta)=(14,6,0)$}
\label{sec:14,6,0}

In this case,
\begin{equation}
S=S_{14,6,0}=U\oplus D_4\oplus D_4\oplus D_4.
\label{latS:14.6.0}
\end{equation}

The lattice $S$ has the basis $\{e,d,f_1,f_2,f_3,f_4,g_1,g_2,g_3,g_4,h_1,h_2,h_3,h_4\}$
where $\{e,d\}$ is the standard basis of $U$,  $f_1,\dots,f_4$, $g_1,\cdots g_4$ and
$h_1,\cdots, h_4$ give
the standard bases of $D_4$.
See \cite{Bou}.

We consider K3 surfaces $X$ with the Picard lattice $S$, the elliptic fibration $|e|$ with the section
$D$ which is a non-singular irreducible rational curve with the class $d$, three degenerate
fibers of the type $\widetilde{\ddd}_4$ with the divisors $F_0+F_1+2F_2+F_3+F_4$, $G_0+G_1+2G_2+G_3+G_4$,
$H_0+H_1+2H_2+H_3+H_4$
where non-singular irreducible rational curves $F_i$, $G_i$, $H_i$, $i=1,...,4$, have classes $f_i$, $g_i$, $h_i$  respectively and $F_0$, $G_0$, $H_0$ have  classes $f_0=e-f_1-2f_2-f_3-f_4$, $g_0=e-g_1-2g_2-g_3-g_4$,
$h_0=e-h_1-2h_2-h_3-h_4$ respectively. They are also non-singular irreducible rational curves.

By considering the canonical involution $\theta$ on $X$, like for the previous cases, we obtain that
the curves $F_2,\ G_2,\ H_2,\ D$  with the classes
\begin{equation}
f_2,\ g_2,\ h_2,\ d
\label{+:14.6.0}
\end{equation}
respectively have the type $+$. Thus they give the curves $E_1,\dots E_k$ where $k=(r-a)/2=4$ in \eqref{Xtheta1}.
 All other non-singular irreducible rational curves on $X$ have the type $-$. Let $C=C^{(1)}$ in
 \eqref{Xtheta1}. Like for the previous cases,
 using this information, we obtain that the class $c$ of $C$ is
 equal to
$$
c=3e^\ast+f_1^\ast+f_3^\ast+f_4^\ast+g_1^\ast+g_3^\ast+g_4^\ast+h_1^\ast+h_2^\ast+h_4^\ast=
$$
\begin{equation}
6e+3d-2f_1-3f_2-2f_3-2f_4-2g_1-3g_2-2g_3-2g_4-2h_1-3h_2-2h_3-2h_4
\label{c:14.6.0}
\end{equation}
where we denote by
$\{e^\ast,d^\ast,f_1^\ast,\dots,f_4^\ast, g_1^\ast,\dots, g_4^\ast, h_1^\ast,\dots, h_4^\ast \}$
the corresponding dual basis from $S^\ast$ for our  basis in the lattice $S$. One can check that $c^2=0$.

We see that  non-singular irreducible rational curves $D$, $F_0$, $F_2$, $G_0$, $G_2$, $H_0$, $H_2$
(or the corresponding classes $d$, $f_0$, $f_2$, $g_0$, $g_2$, $h_0$, $h_2$)
have the Dynkin diagram of the type $\widetilde{\eee_6}$ where $D$ is the central element.
Moreover, $3d+2f_0+f_2+2g_0+g_2+2h_0+h_2=c$ in \eqref{c:14.6.0} where $c$ is the class of
$C=C^{(1)}$ in \eqref{Xtheta1}. Thus, the elliptic pencil $|c|$ has the non-singular fiber $C=C^{(1)}$
in \eqref{Xtheta1} and the degenerate fiber $3D+2F_0+F_2+2G_0+G_2+2H_0+H_2$ of the type $\widetilde{\eee_6}$
with components with classes
\begin{equation}
d,f_0,f_2,g_0,g_2,h_0,h_2
\label{Ecal:14.6.0}
\end{equation}
which contains the curves $E_1,\dots E_k$, $k=4$, in \eqref{Xtheta1}. The elliptic pencil $|c|$ has
the section $F_1$ because $f_1\cdot c=1$ and it has the type $-$. Since the canonical
involution $\theta$ has only two fixed points in $F_1$ corresponding two fibers $C^{(1)}$ and the
degenerate fiber $3D+2F_0+F_2+2G_0+G_2+2H_0+H_2$, the elliptic fibration $|c|$ does not have other
reducible fibers.

Thus, from the theory of elliptic pencils on K3 we obtain that the subgroup of finite index of the automorphism group of
the elliptic fibration $|c|$ is a subgroup of finite index of the group
$MW^\ast$ where
$$
MW=(c,f_1,f_0,g_2,g_0,h_2,h_0,d)^\perp_S.
$$
which has the rank $14-2-6=6$. The lattice $MW$ has no elements with square $(-2)$.

Thus, we obtain

\begin{theorem}
The automorphism group $\Aut X$ of a K3 surface $X$ over $\bc$ with the 2-elementary
Picard lattice $S$ with the invariants
$(r=14,a=6,\delta=0)$, equivalently,
$S=U\oplus 3D_4$,
is equal to the automorphism group of the elliptic fibration $|c|$ of $X$ where $c$ is given in
\eqref{c:14.6.0}. Classes of the curves $E_1,\dots, E_k$ from \eqref{Xtheta1}
are given in \eqref{+:14.6.0} where $k=4$.
The elliptic fibration $|c|$ has
the degenerate fiber of the type $\widetilde{\eee}_{6}$
with components with classes \eqref{Ecal:14.6.0}
and no other reducible fibers.
Then $\Aut X$
is isomorphic to $\bz^6$ up to finite index.

See details in the considerations above.
\label{theorem14.6.0}
\end{theorem}


V.V. Nikulin
\par Steklov Mathematical Institute,
\par ul. Gubkina 8, Moscow 117966, GSP-1, Russia;

\vskip5pt

\par Department of Mathematical Sciences,
\par University of Liverpool,
\par Liverpool, L69 3BX, UK

\vskip5pt

\par E-mail: nikulin@mi-ras.ru\ \ \ vnikulin@liv.ac.uk\ \
vvnikulin@list.ru

\end{document}